# ON VANISHING SUMS OF ROOTS OF UNITY

T. Y. LAM AND K. H. LEUNG

Abstract. Consider the $m$-th roots of unity in $\mathbb{C}$, where $m > 0$ is an integer. We address the following question: For what values of $n$ can one find $n$ such $m$-th roots of unity (with repetitions allowed) adding up to zero? We prove that the answer is exactly the set of linear combinations with non-negative integer coefficients of the prime factors of $m$.

## 1. Introduction

For a given natural number $m$, consider the $m$th roots of unity in the field of complex numbers, $\mathbb{C}$. *For what natural numbers $n$ do there exist $m$th roots of unity $\alpha_1, \cdots, \alpha_n \in \mathbb{C}$ such that $\alpha_1 + \alpha_2 + \cdots + \alpha_n = 0$?* (Such an equation is said to be a vanishing sum of $m$th roots of unity *of weight $n$*.) Although linear relations among roots of unity have been studied rather extensively, a satisfactory answer to the above question is apparently unknown. A couple of explicit examples will show what the set of possible $n$'s can look like. For instance, for $m = 13$, the set of $n$'s is $\{0, 13, 26, 39, \cdots\}$; for $m = 14$, the set of $n$'s is $\{0, 2, 4, 6, 7, 8, 9, 10, \cdots\}$; for $m = 15$, the set of $n$'s is $\{0, 3, 5, 6, 8, 9, 10, 11, \cdots\}$.

For a given $m$, let $W(m)$ be the set of weights $n$ for which there exists a vanishing sum $\alpha_1 + \alpha_2 + \cdots + \alpha_n = 0$, where each $\alpha_i$ is an $m$th root of unity. Since an empty sum is defined to be zero, we agree that $0 \in W(m)$ as in the last paragraph. If $m$ has prime factorization $p_1^{a_1} \cdots p_r^{a_r}$ $(a_i > 0)$, we can easily name a rather large subset of $W(m)$. Indeed, if $\zeta$ is a primitive $p_i$th root of unity (for any $i$), then we have a vanishing sum $1 + \zeta + \cdots + \zeta^{p_i - 1} = 0$ of weight $p_i$. This shows that $W(m)$ contains each $p_i$, and therefore it also contains any linear combination of $p_1, \cdots, p_r$ with non-negative integer coefficients. In particular, if $r \geq 2$, we see that all *sufficiently large* integers $n$ belong to $W(m)$ (for a given $m$). And, if $6 \mid m$, then $W(m)$ consists of all non-negative integers $\neq 1$.

Lam was supported in part by NSF. Research at MSRI is supported in part by NSF grant DMS-9022140.

Leung's research was carried out while he was on sabbatical leave at U.C. Berkeley from the National University of Singapore. The hospitality of the former institution is gratefully acknowledged.

The principal result of this paper is the following:

**Main Theorem.** *For any* $m = p_1^{a_1} \cdots p_r^{a_r}$ *as above, the weight set* $W(m)$ *is exactly given by* $\mathbb{N}\, p_1 + \cdots + \mathbb{N}\, p_r$.

Throughout this paper, $\mathbb{N}$ denotes the semi-ring of non-negative integers. This is a slight deviation from the usual convention that $\mathbb{N} = \{1, 2, 3, \cdots\}$, but it will be convenient for the purposes of this paper. Readers who have misgivings about $0 \in \mathbb{N}$ should feel free to replace our $\mathbb{N}$ by the possibly more reasonable (but obviously clumsy) notation $\mathbb{Z}^+$.

Note that the theorem above implies that $W(m)$ depends only on the prime divisors of $m$, and not on the multiplicities to which they occur in the factorization of $m$. The theorem also shows that any (nonempty) vanishing sum of $m$th roots of unity must have weight $\geq p_1$, where $p_1$ is the smallest prime divisor of $m$. (Vanishing sums of weight $p_1$ turn out to be of the expected type.)

The key technique used for the proof of the Main Theorem is that of group rings. Group rings provide a very natural setting for studying linear relations among roots of unity, but surprisingly they have not been exploited as fully as they should in the literature on the subject. In fact, this is possibly one of the reasons why the result mentioned above has not been discovered earlier. Many of the arguments in this paper would have been rather unwieldy if we were to work with roots of unity alone without the benefit of group rings.

Let $G = \langle z \rangle$ be a cyclic group of order $m$, and let $\zeta$ be a (fixed) primitive $m$th root of unity. There exists a natural ring homomorphism $\varphi$ from the integral group ring $\mathbb{Z}G$ to the ring of cyclotomic integers $\mathbb{Z}[\zeta]$, given by the equation $\varphi(z) = \zeta$. An element of $\mathbb{Z}\,G$, say $x = \sum_{g \in G} x_g\, g$, lies in $\ker(\varphi)$ if and only if $\sum_{g \in G} x_g\, \varphi(g) = 0$ in $\mathbb{Z}[\zeta]$. Therefore, the elements of the ideal $\ker(\varphi)$ correspond precisely to all $\mathbb{Z}$-linear relations among the $m$th roots of unity. For vanishing sums of $m$th roots of unity, we have to look at elements $x = \sum_{g \in G} x_g\, g$ in $\ker(\varphi)$ with all $x_g \geq 0$. In other words, we have to look at $\mathbb{N}G \cap \ker(\varphi)$, where $\mathbb{N}G$ denotes the group semi-ring of $G$ over $\mathbb{N}$. If $x \in \mathbb{N}G \cap \ker(\varphi)$, the weight of the corresponding vanishing sum of $m$th roots of unity is exactly the *augmentation* of the group ring element $x \in \mathbb{Z}G$.

The map $\varphi : \mathbb{Z}G \to \mathbb{Z}[\zeta]$ above will play a central role in this paper, and will often be referred to in the text as "the usual map." It is easy to see that the kernel of $\varphi$ is just the principal ideal generated by $\Phi_m(z)$ in $\mathbb{Z}G$, where $\Phi_m$ is the $m$th cyclotomic polynomial. However, there are other useful descriptions of $\ker(\varphi)$. For instance, a theorem of Rédei [$\mathrm{R}_1$, $\mathrm{R}_2$], de Bruijn [deB] and Schoenberg [Sch] can be recast in the language of group rings to give a natural family of ideal generators for $\ker(\varphi)$ in terms of the minimal subgroups of $G$. In §2, we give a new proof of this theorem using induction and group-theoretic techniques. In §3, we use similar methods to prove several other useful facts about $\mathbb{N}G \cap \ker(\varphi)$. After this preparatory work, we prove in §4 a Lower Bound Theorem (4.8) on the weights of minimal vanishing sums. From this result, the Main Theorem (5.2) follows easily. In §6, we prove the existence and uniqueness of minimal vanishing sums with the smallest weight predicted by the Lower Bound Theorem. The final section §7 offers an application of the Main Theorem to the character theory of finite groups. Throughout this paper, we work in characteristic 0. The study of the same problems in characteristic $p$ requires different techniques, and will be reported elsewhere (see [LL]).

Aside from being of intrinsic interest, vanishing sums of roots of unity arise naturally also in a number of algebraic, geometric and combinatorial contexts: for instance, in cyclotomy and difference sets [St], factorization problems in groups [deB], trigonometric diophantine equations [CJ], and in the study of polar rational polygons (convex polygons with integral sides whose angles are rational when measured in degrees) [Sch], [Ma]. For a partial survey of the literature up to 1978, see [Le]. For obvious reasons, it is of interest to study vanishing sums $\alpha_1 + \cdots + \alpha_n = 0$ which are *minimal*, in the sense that no proper subsums thereof can be zero. Minimal vanishing sums involving "few" distinct roots of unity were classified by Mann [Ma] and Conway-Jones [CJ]. Recently, in connection with their work on counting the intersection points of the diagonals of a regular polygon, Poonen and Rubinstein [PR] have classified all minimal vanishing sums $\alpha_1 + \cdots + \alpha_n = 0$ of weight $n \leq 12$. Of course, we can always multiply a vanishing sum by a root of unity to get another; we say that the latter is *similar* to the former, or that it is obtained from the former "by a rotation". Naturally, the classification of minimal vanishing sums needs to be done only *up to rotations* (by roots of unity).

To conclude this Introduction, a few notational remarks are in order. For any natural number $m$, $\zeta_m$ shall denote a *primitive* $m$th root of unity (in $\mathbb{C}$). For any commutative ring (resp. semi-ring) $k$ and any group $G$, we shall write $kG$, or sometimes $k[G]$, for the group ring (resp. group semi-ring) of $G$ over $k$. Elements of $kG$ will be written in the form $x = \sum_{g \in G} x_g\, g$, where $x_g \in k$ are zero except for a finite number of $g$'s. The number of nonzero coefficients $x_g$ is denoted by $\varepsilon_0(x)$, and the sum of such coefficients is denoted by $\varepsilon(x)$. The latter is called the *augmentation* of $x$, and $x \mapsto \varepsilon(x)$ defines a $k$-homomorphism $\varepsilon : kG \to k$, called the augmentation map. For any *finite* subset $H \subseteq G$, we shall write $\sigma(H)$ for the sum $\sum_{h \in H} h$ in the group (semi-)ring $kG$. Two basic properties of $\sigma(H)$ to be used freely in the sequel are that $\varepsilon(\sigma(H)) = \varepsilon_0(\sigma(H)) = |H|$ (the cardinality of $H$), and that, if $H$ is a subgroup, $\sigma(H) \cdot h = \sigma(H)$ for any $h \in H$. In the case when $k = \mathbb{Z}$,

we can define a partial ordering on $\mathbb{Z}G$, by declaring that

$$y = \sum y_g\, g \geq x = \sum x_g\, g \iff y_g \geq x_g \quad \text{for every } g \in G. \tag{1.1}$$

Note that $y \geq x$ iff $y - x \geq 0$, and the "positive cone" $\{z \in \mathbb{Z}G: \ z \geq 0\}$ is precisely the group semi-ring $\mathbb{N}G$.

*Acknowledgment.* We thank R. Guralnick and H. Lenstra for helpful comments on this work, and B. Poonen and M. Rubinstein for providing their preprint [PR].

## 2. The Rédei-de Bruijn-Schoenberg Theorem

In this section, we recast a theorem of Rédei, de Bruijn and Schoenberg in the setting of group rings, and give an easy inductive proof for this theorem. This work is in preparation for what is to come in §§3-5.

For a cyclic group $G = \langle z \rangle$ of order $m = p_1^{a_1} \cdots p_r^{a_r}$ and a primitive $m$th root of unity $\zeta = \zeta_m$, consider the usual map

$$\varphi: \ \mathbb{Z}G \to \mathbb{Z}[\zeta], \quad \text{with } \varphi(z) = \zeta. \tag{2.1}$$

As we have observed before, $\ker(\varphi) = \mathbb{Z}G \cdot \Phi_m(z)$, where $\Phi_m$ denotes the $m$th cyclotomic polynomial, i.e. the minimal polynomial of $\zeta$ over $\mathbb{Q}$. It turns out that there is another useful description of $\ker(\varphi)$. Let $P_i$ $(1 \leq i \leq r)$ be the unique subgroup of order $p_i$ in $G$, and consider $\sigma(P_i) := \sum_{g \in P_i} g \in \mathbb{Z}G$. For any nonidentity element $g \in P_i$, we have $\sigma(P_i) \cdot g = \sigma(P_i)$, so we must have $\sigma(P_i) \in \ker(\varphi)$, since $\varphi(g) \neq 1$ and $\mathbb{Z}[\zeta]$ is an integral domain. The following result calculates $\ker(\varphi)$ in terms of the special elements $\sigma(P_1), \cdots, \sigma(P_r)$.

**Theorem 2.2.** (*cf.* [$R_1$: Hilfssatz 4][1], [deB: Th. 1], [Sch: Th. 1])

$$\ker(\varphi) = \sum_{i=1}^{r} \mathbb{Z}G \cdot \sigma(P_i), \quad \textit{and} \quad \ker(\varphi) = \mathbb{Z} \cdot \sigma(P_1) \ \textit{in case} \ r = 1.$$

**Proof.** We proceed by induction on $r$. First assume $r = 1$ and write $p = p_1$, $a = a_1$. Then, $m = p^a$, $P_1 = \langle z^{p^{a-1}} \rangle$, and we have

$$\sigma(P_1) = 1 + z^{p^{a-1}} + (z^{p^{a-1}})^2 + \cdots + (z^{p^{a-1}})^{p-1}.$$

On the other hand,

$$\begin{aligned} \Phi_m(X) &= \left(X^{p^a} - 1\right) / \left(X^{p^{a-1}} - 1\right) \\ &= \left(X^{p^{a-1}}\right)^{p-1} + \left(X^{p^{a-1}}\right)^{p-2} + \cdots + X^{p^{a-1}} + 1. \end{aligned}$$

[1] As pointed out by de Bruijn, Rédei's proof of his Hilfssatz 4 in [$R_1$] was incomplete. Complete proofs appeared later in [deB] and [$R_2$].

Therefore, $\Phi_m(z)$ is exactly $\sigma(P_1)$, and we have

$$\ker(\varphi) = \mathbb{Z}G \cdot \Phi_m(z) = \mathbb{Z}G \cdot \sigma(P_1) = \mathbb{Z} \cdot \sigma(P_1).$$

For $r \geq 2$, write $G = H \times H'$ where $n := |H| > 1$, $n' := |H'| > 1$, and $(n, n') = 1$. As in (2.1), we have surjections $\psi : \mathbb{Z}H \to \mathbb{Z}[\zeta_n]$ and $\psi' : \mathbb{Z}H' \to \mathbb{Z}[\zeta_{n'}]$ (where we may asume $\zeta_n \zeta_{n'} = \zeta_m$). Let $I = \ker(\psi)$ and $I' = \ker(\psi')$. Since $\mathbb{Z}[\zeta_n]$ is $\mathbb{Z}$-free, there exists a $\mathbb{Z}$-basis $\{e_i, f_j\}$ for $\mathbb{Z}H$ such that $\{e_i\}$ is a $\mathbb{Z}$-basis for $I$ and $\{\psi(f_j)\}$ is a $\mathbb{Z}$-basis for $\mathbb{Z}[\zeta_n]$. Similarly, we fix a basis $\{e'_k, f'_\ell\}$ for $\mathbb{Z}H'$. Then, $\mathbb{Z}H \otimes \mathbb{Z}H'$ has a $\mathbb{Z}$-basis

$$\{e_i \otimes e'_k,\ e_i \otimes f'_\ell,\ f_j \otimes e'_k,\ f_j \otimes f'_\ell\},$$

where the first three sets of elements lie in $\ker(\psi \otimes \psi')$. Since $\{\psi(f_j) \otimes \psi'(f'_\ell)\}$ is a $\mathbb{Z}$-basis for $\mathbb{Z}[\zeta_n] \otimes \mathbb{Z}[\zeta_{n'}]$, it follows that $\{e_i \otimes e'_k, e_i \otimes f'_\ell, f_j \otimes e'_k\}$ is a $\mathbb{Z}$-basis for $\ker(\psi \otimes \psi')$. Now, identifying $\mathbb{Z}H \otimes \mathbb{Z}H'$ with $\mathbb{Z}G$ and $\mathbb{Z}[\zeta_n] \otimes \mathbb{Z}[\zeta_{n'}]$ with $\mathbb{Z}[\zeta_m]$, we see that

$$\ker(\varphi) = \ker(\psi \otimes \psi') = I \otimes \mathbb{Z}H' + \mathbb{Z}H \otimes I'. \tag{2.3}$$

From this, the Theorem follows immediately by induction. □

In terms of roots of unity, (2.2) says that any $\mathbb{Z}$-linear relation among the $m$th roots of unity can be obtained from the basic relations

$$1 + \zeta_{p_i} + \cdots + \zeta_{p_i}^{p_i - 1} = 0 \quad (1 \leq i \leq r) \tag{2.4}$$

by addition, subtraction, and rotation. However, this does not mean that every *vanishing sum* of $m$th roots of unity can be obtained from those of the type (2.4) by addition and rotation. In other words, there exist in general minimal vanishing sums which are *not* similar to those in (2.4).

**Example 2.5.** Let $m$ be an integer with at least three prime factors $p$, $q$, $\ell$, and let $\alpha = \zeta_p$, $\beta = \zeta_q$ and $\gamma = \zeta_\ell$. Following a construction of Rédei [Re$_3$: Satz 9], consider the sum

$$(\alpha + \cdots + \alpha^{p-1})(\beta + \cdots + \beta^{q-1}) + \gamma + \cdots + \gamma^{\ell-1} = (-1)(-1) + (-1) = 0. \tag{2.6}$$

We claim that this vanishing sum is minimal. To see this, consider any vanishing subsum, say

$$\alpha\, b_1 + \alpha^2\, b_2 + \cdots + \alpha^{p-1}\, b_{p-1} + c = 0, \tag{2.7}$$

where each $b_i$ is a subsum of $\beta + \cdots + \beta^{q-1}$, and $c$ is a subsum of $\gamma + \cdots + \gamma^{\ell-1}$. If $c$ is the empty sum, then (2.7) and the linear disjointness of $\mathbb{Q}(\alpha)$ and $\mathbb{Q}(\beta)$ over $\mathbb{Q}$ show that all $b_i = 0$. This implies that each $b_i$ (as well as $c$) is an empty sum, so we are done. If $c$ is not the empty sum, then the linear disjointness of $\mathbb{Q}(\gamma)$ from $\mathbb{Q}(\alpha, \beta)$ implies that $c$ is the entire sum $\gamma + \cdots + \gamma^{\ell-1} = -1$, and (2.7) and the linear disjointness of $\mathbb{Q}(\alpha)$ from $\mathbb{Q}(\beta)$ imply in turn that all $b_i = -1$. This is possible only if each $b_i$ is also the entire sum $\beta + \cdots + \beta^{q-1}$, as desired. Clearly, the minimal sum (2.6) is *not* similar to one of the type (2.4). Note that, in the special case $p = 2$, we have $\alpha = -1$; here (2.6) takes on the simpler form

$$-\beta - \cdots - \beta^{q-1} + \gamma_1 + \cdots + \gamma^{\ell-1} = 0. \tag{2.8}$$

Minimal sums of this type were used recently by Poonen and Rubinstein in their study of the intersection points of the diagonals of a regular polygon [PR: §3].

Vanishing sums of the form (2.6) turn out to have a special significance. As we shall see later in (6.5), they have the smallest possible weight among all "asymmetric" minimal vanishing sums involving $m$th roots of unity, for any $m$ whose smallest prime divisors are $p$, $q$, $\ell$.

## 3. Results on $\mathbb{N}G \cap \ker(\varphi)$

In this section, we study the map $\varphi$ in (2.1) for a finite cyclic group $G$, and prove a few results on the structure of elements lying in the intersection $\mathbb{N}G \cap \ker(\varphi)$. The first result below is essentially equivalent to [CJ: Theorem 1], but our proof is different from that of Conway and Jones, and the result is stated here in terms of group rings, in a form most convenient for later use.

**Theorem 3.1.** *Let $G$ be a cyclic group of order $m = p_1^{a_1} \cdots p_r^{a_r}$ where $p_1, \cdots, p_r$ are distinct primes, and let $\varphi : \mathbb{Z}G \to \mathbb{Z}[\zeta]$ be the usual map, where $\zeta = \zeta_m$. Let $G_0 \subseteq G$ be the (unique) subgroup of order $p_1 \cdots p_r$, and let $\{g_j : 1 \le j \le [G : G_0]\}$ be a complete set of coset representatives of $G$ with respect to $G_0$. Then $\mathbb{N}\,G \cap \ker(\varphi) = \sum_j g_j\,(\mathbb{N}G_0 \cap \ker(\varphi))$.*

**Proof.** We need only prove the inclusion "$\subseteq$". Let $x \in \mathbb{N}G \cap \ker(\varphi)$. By (2.2), we can write $x = \sum x_i\,\sigma(P_i)$, where $P_i$ is the unique subgroup of order $p_i$ in $G_0$, and $x_i \in \mathbb{Z}G$. Using the decomposition $\mathbb{Z}G = \sum_j g_j\,\mathbb{Z}G_0$, we can write $x_i = \sum_j g_j\,y_{ij}$ where $y_{ij} \in \mathbb{Z}G_0$. Then

$$x = \sum_i x_i\,\sigma(P_i) = \sum_i \sum_j g_j y_{ij}\,\sigma(P_i) = \sum_j g_j z_j,$$

where $z_j = \sum_i y_{ij}\,\sigma(P_i) \in \mathbb{Z}G_0$. Since the sum $\sum_j g_j\,\mathbb{Z}G_0 = \mathbb{Z}G$ is direct, the fact that $x \in \mathbb{N}G$ implies that $z_j \in \mathbb{N}G_0$ for all $j$. On the other hand, $\varphi(z_j) = \sum_i \varphi(y_{ij})\varphi(\sigma(P_i)) = 0$, so $z_j \in \ker(\varphi)$ for all $j$. Therefore, we have $x \in \sum_j g_j\,(\mathbb{N}\,G_0 \cap \ker(\varphi))$, as desired. $\square$

**Corollary 3.2.** *If $\alpha_1 + \cdots + \alpha_n = 0$ is a minimal vanishing sum of $m$ th roots of unity, then after a suitable rotation, we may assume that all $\alpha_i$ 's are $m_0$ th roots of unity where $m_0$ is square-free.*

**Proof.** Suppose the given relation $\alpha_1 + \cdots + \alpha_n = 0$ corresponds to an $x \in \mathbb{N}G \cap \ker(\varphi)$, where $G$ is a group as in (3.1). Using the notations there, we have a decomposition $x = \sum_j g_j z_j$ where $z_j \in \mathbb{N}G_0 \cap \ker(\varphi)$. Since the given relation is minimal, we must have $x = g_j z_j$ for some $j$. Therefore, after a rotation by $\varphi(g_j)^{-1}$, the given relation becomes $\alpha'_1 + \cdots + \alpha'_n = 0$ where the $\alpha'_i$'s are $m_0$ th roots of unity with $m_0 = |G_0|$ square-free. □

In the case when $m$ has at most two prime divisors, there is a very explicit description of $\mathbb{N}G \cap \ker(\varphi)$. This result can be traced back to the work of de Bruijn [DeB: §3]. For the sake of completeness, we offer here a self-contained proof. In fact, this proof in terms of group rings (for the case $r = 2$) will set the stage for several of the inductive proofs to be given in the next section.

**Theorem 3.3.** *Keep the notations of Theorem 3.1, and let $P_i$ denote the unique subgroup of order $p_i$ in $G$. (1) If $r = 1$, $\mathbb{N}G \cap \ker(\varphi) = \mathbb{N} \cdot \sigma(P_1)$. (2) If $r = 2$, $\mathbb{N}G \cap \ker(\varphi) = \mathbb{N}P_1 \cdot \sigma(P_2) + \mathbb{N}P_2 \cdot \sigma(P_1)$.*

**Proof.** (1) follows from (2.2) (in the case $r = 1$). Now assume $r = 2$. In view of Theorem 3.1, it suffices to prove (2) for the unique subgroup $G_0$ of order $p_1 p_2$ in $G$. Let us assume, therefore, that $|G| = |G_0| = qp$, where $q = p_1$, $p = p_2$. Let $P_2 = \langle g \rangle$, so that $G = P_1 \times \langle g \rangle$. For any $x \in \mathbb{N}G \cap \ker(\varphi)$, we can write $x = x_0 + x_1 g + \cdots + x_{p-1} g^{p-1}$, where $x_i \in \mathbb{N}P_1$. Then, $x \in \ker(\varphi)$ implies that

$$\varphi(x_0) + \varphi(x_1)\,\zeta_p + \cdots + \varphi(x_{p-1})\,\zeta_p^{p-1} = 0,$$

where $\zeta_p$ is a primitive $p$ th root of unity. Since $\mathbb{Q}(\zeta_p)$ is linearly disjoint from $\mathbb{Q}(\zeta_q)$ over $\mathbb{Q}$ (by a theorem of Kronecker), we must have $\varphi(x_0) = \varphi(x_1) = \cdots = \varphi(x_{p-1})$ in $\mathbb{Q}(\zeta_q)$. Say $\varepsilon(x_i)$ is the smallest among all $\varepsilon(x_j)$'s. From $\varphi(x_j - x_i) = 0$, we have $x_j - x_i = z_j\,\sigma(P_1)$ for some $z_j \in \mathbb{Z}$ (by (2.1) for one prime). Then

$$0 \le \varepsilon(x_j) - \varepsilon(x_i) = z_j \cdot \varepsilon(\sigma(P_1)) = z_j\, q$$

implies that each $z_j \ge 0$. Therefore,

$$\begin{aligned} x &= x_0 + x_1 g + \cdots + x_{p-1} g^{p-1} \\ &= (x_i + z_0\,\sigma(P_1)) + (x_i + z_1\,\sigma(P_1))g + \cdots + (x_i + z_{p-1}\,\sigma(P_1))g^{p-1} \\ &= x_i(1 + g + \cdots + g^{p-1}) + (z_0 + z_1 g + \cdots + z_{p-1} g^{p-1})\,\sigma(P_1), \end{aligned}$$

which lies in $\mathbb{N}P_1 \cdot \sigma(P_2) + \mathbb{N}P_2 \cdot \sigma(P_1)$, as desired. □

The following is a direct consequence of (3.3)(2). In the special case when $m = 2p$, this was noted recently by Poonen and Rubinstein in [PR: Lemma 2].

**Corollary 3.4.** *Let $m = p^a q^b$, where $p, q$ are primes. Then, up to a rotation, the only minimal vanishing sums of $m$th roots of unity are: $1 + \zeta_p + \cdots + \zeta_p^{p-1} = 0$, and $1 + \zeta_q + \cdots + \zeta_q^{q-1} = 0$.*

## 4. The Lower Bound Theorem

In this and the following sections, we'll keep the notations set up in the Introduction and at the beginning of §2. We say that a nonzero element $x \in \mathbb{N}G \cap \ker(\varphi)$ is *minimal* if it cannot be decomposed into a sum of two nonzero elements in $\mathbb{N}G \cap \ker(\varphi)$. In other words, $x$ is minimal if and only if $\varphi(x) = 0$ represents a minimal vanishing sum of $m$th roots of unity. For $1 \le i \le r$, let $P_i$ be the unique subgroup of order $p_i$ in $G$. Then, for any $g \in G$, $g \cdot \sigma(P_i) \in \mathbb{N}G \cap \ker(\varphi)$ is minimal, since it corresponds (up to a rotation) to the minimal vanishing sum

$$1 + \zeta_{p_i} + \cdots + \zeta_{p_i}^{p_i - 1} = 0.$$

We shall refer to $\{g \cdot \sigma(P_i)\}$ as the *symmetric* minimal elements of $\mathbb{N}G \cap \ker(\varphi)$; the other minimal elements will be referred to as the *asymmetric* ones. A similar terminology will be used for minimal vanishing sums.

In the case $r \le 2$, (3.3) implies that all minimal elements in $\mathbb{N}G \cap \ker(\varphi)$ are symmetric. However, when $r \ge 3$, (2.5) shows that there exist asymmetric minimal elements. In this section, we shall study these elements in $\mathbb{N}G \cap \ker(\varphi)$, where $G$ is a cyclic group of order $m = p_1^{a_1} \cdots p_r^{a_r}$ $(r \ge 3)$. The main result here is Theorem 4.8 which provides an effective lower bound on the $\varepsilon_0$ (size of the support) of such asymmetric (minimal) elements. To begin with, we prove a preliminary result (in the case when $|G|$ is square-free) on the $\varepsilon_0$ of two elements $x, y \in \mathbb{N}G$ which have the same image under the homomorphism $\varphi$. Recall that a partial ordering "$\ge$" for elements in $\mathbb{Z}G$ was defined in (1.1).

**Theorem 4.1.** *Let $G$ be a cyclic group of order $m = p_1 p_2 \cdots p_r$ where $p_1 < p_2 < \cdots < p_r$ are primes and $r \ge 2$. Let $\varphi : \mathbb{Z}G \to \mathbb{Z}[\zeta]$ be the usual map, where $\zeta = \zeta_m$. Let $x, y \in \mathbb{N}G$ be such that $\varphi(x) = \varphi(y)$. If $\varepsilon_0(x) \le p_1 - 1$, then we have either* (A) *$y \ge x$, or* (B) *$\varepsilon_0(y) \ge (p_1 - \varepsilon_0(x))(p_2 - 1)$. In Case* (A), *we have $\varepsilon_0(y) \ge \varepsilon_0(x)$, and in Case* (B), we have $\varepsilon_0(y) > \varepsilon_0(x)$.

To better understand this Theorem, a simple illustrative example is in order.

**Example 4.2.** Let $|G| = p_1 p_2$, where $p_1 < p_2$ are primes. Let $P_i \subseteq G$ be the unique subgroup of order $p_i$, and $P_i^* = P_i \backslash \{1\}$. Let $P_1 = X \cup X'$ be any partition of $P_1$, with $X, X' \ne \emptyset$. Now let $x = c \cdot \sigma(X)$, and $y = c \cdot \sigma(X')\,\sigma(P_2^*)$, where $c$ is any positive integer. Then $\varepsilon_0(x) = |X| \le p_1 - 1$, and since $\varphi(\sigma(X)) + \varphi(\sigma(X')) = \varphi(\sigma(P_1)) = 0$, we have

$$\varphi(y) = c \cdot \varphi(\sigma(X')) \cdot \varphi(\sigma(P_2^*)) = -c \cdot \varphi(\sigma(X))(-1) = \varphi(x),$$

checking the hypotheses in the Theorem. In this example, (A) clearly does not hold, and (B) holds *with an equality*, since

$$\varepsilon_0(y) = \varepsilon_0(\sigma(X'))\, \varepsilon_0(\sigma(P_2^*)) = (p_1 - |X|) \cdot |P_2^*| = (p_1 - \varepsilon_0(x))(p_2 - 1).$$

This example shows that in general the conclusion in (4.1) is the best possible. Note that in the special case when $X = \{1\}$, we have $x = c$ and $y = c \cdot \sigma(P_1^*)\, \sigma(P_2^*)$.

**Proof of (4.1).** The last statement in the theorem follows since, in Case (B), we'll have

$$\varepsilon_0(y) \ge (p_1 - \varepsilon_0(x))(p_2 - 1) \ge p_2 - 1 > p_1 - 1 \ge \varepsilon_0(x).$$

The proof of the theorem will be by induction on $r \ge 2$. Let $H \subseteq G$ be the (unique) subgroup of order $p_1 \cdots p_{r-1}$, and let $g \in G$ be an element of order $p := p_r$, so that $G = H \times \langle g \rangle$. Then there are unique expressions

$$\begin{aligned} x &= x_0 + x_1 g + \cdots + x_{p-1} g^{p-1}, \\ y &= y_0 + y_1 g + \cdots + y_{p-1} g^{p-1}, \end{aligned}$$

where $x_i,\, y_i \in \mathbb{N}H$. Let $I = \{i : x_i = 0\}$. This is a nonempty set since $\varepsilon_0(x) \le p_1 - 1 < p - 1$. In the set $\{y_i : i \in I\}$, choose $y_j$ such that $\varepsilon_0(y_j)$ is the smallest. From the hypothesis $\varphi(x) = \varphi(y)$, we have $\sum_{i=0}^{p-1} \varphi(x_i - y_i)\zeta_p^i = 0$, where, as usual, $\zeta_p$ denotes a primitive $p$ th root of unity. Since $\varphi(x_i - y_i) \in \mathbb{Q}(\zeta_{|H|})$, and $\mathbb{Q}(\zeta_{|H|})$, $\mathbb{Q}(\zeta_p)$ are linearly disjoint over $\mathbb{Q}$, we must have $\varphi(y_i - x_i) = \varphi(y_j - x_j)$ for all $i$, or, equivalently,

$$\varphi(y_i) = \varphi(x_i + y_j) \quad \text{for all } i. \tag{4.3}$$

Choose $k$ such that $\varepsilon_0(x_k)$ is maximum (among all $\varepsilon_0(x_i)$ 's). We shall distinguish the following two main cases.

*Case 1.* $\varepsilon_0(x_k) + \varepsilon_0(y_j) \ge p_1$. Let $t := p - |I|$, which is the number of nonzero $x_i$ 's. We may assume that $t \ge 1$, for otherwise $x = 0$ and $y \ge x$ holds. Note the following obvious upper and lower bounds on $\varepsilon_0(x)$:

$$\varepsilon_0(x_k) + t - 1 \le \varepsilon_0(x) \le \varepsilon_0(x_k)\, t.$$

Using the definition of $y_j$, we have

$$\begin{aligned}
\varepsilon_0(y) &\geq |I|\cdot\varepsilon_0(y_j) = (p-t)\,\varepsilon_0(y_j)\\
&\geq (p_2-t)(p_1-\varepsilon_0(x_k))\\
&= p_1p_2 - tp_1 - \varepsilon_0(x_k)\,p_2 + \varepsilon_0(x_k)\,t\\
&= p_1p_2 + t(p_2-p_1) - p_2 - (\varepsilon_0(x_k)+t-1)\,p_2 + \varepsilon_0(x_k)\,t\\
&\geq p_1p_2 + (p_2-p_1) - p_2 - \varepsilon_0(x)\,p_2 + \varepsilon_0(x)\\
&= (p_1-\varepsilon_0(x))\,(p_2-1),
\end{aligned}$$

so we have proved (B) in this case.

*Case 2.* $\varepsilon_0(x_k)+\varepsilon_0(y_j) \leq p_1-1$. This case assumption means that $\varepsilon_0(x_i)+\varepsilon_0(y_j) \leq p_1-1$ for all $i$. We shall first take care of the case $r=2$ (to start the induction). In this case, $|H|=p_1$, so by (4.3) and (the one-prime case of) (2.2),

$$y_i = x_i + y_j + z_i\sigma(H) \quad \text{for some } z_i \in \mathbb{Z}. \tag{4.4}$$

If some $z_i<0$, then $x_i+y_j = y_i+|z_i|\cdot\sigma(H)$ implies that $\varepsilon_0(x_i)+\varepsilon_0(y_j) \geq \varepsilon_0(x_i+y_j) = p_1$, a contradiction. Therefore, we must have $z_i \geq 0$ for all $i$. It follows from (4.4) that $y_i \geq x_i$ for all $i$, and hence $y \geq x$ in this case.

Assume now $r \geq 3$. Since $\varphi(y_i)=\varphi(x_i+y_j)$ and $\varepsilon_0(x_i+y_j) \leq \varepsilon_0(x_i)+\varepsilon_0(y_j) \leq p_1-1$, we can apply the inductive hypothesis to the pair $y_i$ and $x_i+y_j$ in $\mathbb{N}H$. In particular, we will have

$$\varepsilon_0(y_i) \geq \varepsilon_0(x_i+y_j) \quad \text{for all } i. \tag{4.5}$$

If $y_i \geq x_i+y_j$ for all $i$, then $y_i \geq x_i$ for all $i$, and we have $y \geq x$, proving (A) in this case. Otherwise, our inductive hypothesis implies that there exists an $\ell$ such that

$$\varepsilon_0(y_\ell) \geq (p_1-\varepsilon_0(x_\ell+y_j))(p_2-1). \tag{4.6}$$

Note that, from (4.5), $\varepsilon_0(y_i) \geq \varepsilon_0(y_j)$ for all $i$. Using this, we have

$$\begin{aligned}
\varepsilon_0(y) &= \varepsilon_0(y_\ell) + \sum_{i\neq\ell}\varepsilon_0(y_i)\\
&\geq (p_1-\varepsilon_0(x_\ell+y_j))(p_2-1) + (p-1)\,\varepsilon_0(y_j)\\
&\geq (p_1-\varepsilon_0(x_\ell))(p_2-1) + (p-p_2)\,\varepsilon_0(y_j)\\
&\geq (p_1-\varepsilon_0(x))(p_2-1),
\end{aligned}$$

proving (B) in this case. □

**Corollary 4.7.** *Theorem 4.1 holds verbatim with $\varepsilon_0$ replaced throughout by the augmentation $\varepsilon$.*

**Proof.** It is easy to check that every step of the above proof goes through if we use the augmentation $\varepsilon$ instead of $\varepsilon_0$. Alternatively, we may note that (4.1) is actually stronger than (4.7). For, if we assume that $\varepsilon(x) \le p_1 - 1$, then $\varepsilon_0(x) \le p_1 - 1$ also, so we have either (A) or (B) in (4.1). In the case (B), we'll have

$$\varepsilon(y) \ge \varepsilon_0(y) \ge (p_1 - \varepsilon_0(x))(p_2 - 1) \ge (p_1 - \varepsilon(x))(p_2 - 1). \qquad \square$$

We are now ready to establish a lower bound (in terms of the $p_i$'s) for the $\varepsilon_0$ of the asymmetric minimal elements in $\mathbb{N}G \cap \ker(\varphi)$. This crucial result, coupled with a simple fact from elementary number theory, will lead quickly to a proof of the Main Theorem stated in the Introduction.

**Lower Bound Theorem 4.8.** *Let $G$ be a cyclic group of order $m = p_1^{a_1} \cdots p_r^{a_r}$, where $p_1 < \cdots < p_r$ are primes, and let $\mathbb{Z}G \to \mathbb{Z}[\zeta]$ be the usual map, where $\zeta = \zeta_m$. For any minimal element $x \in \mathbb{N}G \cap \ker(\varphi)$, we have either* (A) *$x$ is symmetric, or* (B) *$r \ge 3$ and $\varepsilon(x) \ge \varepsilon_0(x) \ge p_1(p_2 - 1) + p_3 - p_2 > p_3$.*

**Proof.** By Theorem 3.1, we may assume that all the exponents $a_i$ are 1. The proof will be again by induction on $r$. In the case $r \le 2$, (3.3) implies that $x$ is necessarily symmetric, so (A) always holds in this case. This starts the induction, and we may now proceed to the case $r \ge 3$.

Write $x = x_0 + x_1 g + \cdots + x_{p-1} g^{p-1}$ as in the proof of (4.1), where $g$ has order $p := p_r$, $x_k \in \mathbb{N}H$, and $|H| = p_1 \cdots p_{r-1}$. Since $\varphi(x) = 0$, we have $\varphi(x_0) = \varphi(x_1) = \cdots = \varphi(x_{p-1})$ as before (by the linear disjointness argument). Choose $i$ such that $\varepsilon_0(x_i)$ is the smallest. We shall argue in the following three cases.

*Case 1.* $\varepsilon_0(x_i) \ge p_1$. In this case, we have

$$\begin{aligned} \varepsilon_0(x) \;&\ge\; \varepsilon_0(x_i)\, p \;\ge\; p_1 p_3 \;=\; p_1 (p_2 + p_3 - p_2) \\ &>\; p_1 p_2 + p_3 - p_2 \;>\; p_1 (p_2 - 1) + p_3 - p_2. \end{aligned}$$

*Case 2.* $\varepsilon_0(x_i) = 0$. This means that $x_i = 0$, so we have $\varphi(x_k) = \varphi(x_i) = 0$ for all $k$, i.e. $x_k \in \mathbb{N}H \cap \ker(\varphi)$. Since $x = x_0 + x_1 g + \cdots + x_{p-1} g^{p-1}$ is minimal, we must have $x = x_k\, g^k$ for some $k$, with necessarily $x_k$ minimal in $\mathbb{N}H \cap \ker(\varphi)$. Invoking the inductive hypothesis, $x_k$ is either symmetric, or we have $r - 1 \ge 3$ and $\varepsilon_0(x) = \varepsilon_0(x_k) \ge p_1(p_2 - 1) + p_3 - p_2$, as desired.

*Case 3.* We may assume now that $1 \le \varepsilon_0(x_i) \le p_1 - 1$. By (4.1) (applied to the elements $x_i$, $x_j \in \mathbb{N}H$ ), we have the following two possibilities:

*Subcase 1.* $x_j \ge x_i$ for all $j$. In this case,

$$x = x_0 + x_1 g + \cdots + x_{p-1} g^{p-1} \ge x_i + x_i g + \cdots + x_i g^{p-1} = x_i\, \sigma(\langle g \rangle).$$

Since $x$ is minimal, we must have $x = x_i\,\sigma(\langle g\rangle)$ and $x_i \in H$, so $x$ is symmetric in this case.

*Subcase 2.* There exists $j$ such that $\varepsilon_0(x_j) \geq (p_1 - \varepsilon_0(x_i))(p_2 - 1)$. In this case,

$$\begin{aligned} \varepsilon_0(x) &= \varepsilon_0(x_j) + \sum_{k \neq j} \varepsilon_0(x_k) \\ &\geq (p_1 - \varepsilon_0(x_i))(p_2 - 1) + (p-1)\varepsilon_0(x_i) \\ &= p_1(p_2 - 1) + (p - p_2)\varepsilon_0(x_i) \\ &\geq p_1(p_2 - 1) + p - p_2 \\ &\geq p_1(p_2 - 1) + p_3 - p_2. \end{aligned}$$

In any case, we have now shown that either (A) or (B) holds. (For the last inequality in (B), note that $p_1(p_2-1)+p_3-p_2 = p_1p_2-p_2-p_1+p_3 \geq (p_2-p_1)+p_3 > p_3$.) □

**Corollary 4.9.** *In the notations of* (4.8), *any element* $u \in \mathbb{N}G \cap \ker(\varphi)$ *with* $\varepsilon_0(u) < p_1(p_2-1)+p_3-p_2$ *lies in* $\sum_i \mathbb{N}G \cdot \sigma(P_i)$, *where* $P_i$ *is the subgroup of order* $p_i$ *in* $G$.

More can be said about the Lower Bound Theorem (4.8). But at this point, it is perhaps imperative to show first how the Main Theorem can be deduced from it. After showing this in the next section, we shall return in §6 to the Lower Bound Theorem, and determine the structure of the asymmetric minimal elements of the smallest support (resp. weight) in $\mathbb{N}G \cap \ker(\varphi)$.

## 5. The Main Theorem

After all the preparation in the previous sections, it is now an easy matter to prove the Main Theorem. We need just one more elementary number-theoretic fact, the proof of which we shall leave to the reader.

**Lemma 5.1.** (See [LeV: p.22, Ex.4]) *Let* $p$, $q$ *be relatively prime positive integers. If* $n$ *is any integer* $\geq (p-1)(q-1)$, *then* $n \in \mathbb{N}\,p + \mathbb{N}\,q$.

The Main Theorem stated in the Introduction of this paper concerns the computation of the set $W(m)$ of integers $n$ for which there exists (in $\mathbb{C}$) a vanishing sum of $m$th roots of unity of weight $n$. We shall now restate this result in the convenient language of group rings, and derive it as a byproduct of Theorem 4.8.

**Theorem 5.2.** *Keep the notations in* (4.8), *and let* $x \in \mathbb{N}G \cap \ker(\varphi)$. *Then* $\varepsilon(x) \in \sum_{i=1}^{r} \mathbb{N}\,p_i$. *In other words,* $W(m) = \sum_{i=1}^{r} \mathbb{N}\,p_i$.

**Proof.** Since $x$ can be decomposed into a sum of minimal elements in $\mathbb{N}G \cap \ker(\varphi)$, it suffices to prove the theorem for minimal elements $x$. By Theorem 4.8, $x$ is either symmetric, or we'll have $r \geq 3$ and $\varepsilon(x) \geq p_1(p_2-1)+p_3-p_2$. In the former case, $\varepsilon(x) = p_i$ for some $i$. In the latter case,

$$\varepsilon(x) > p_1(p_2-1) > (p_1-1)(p_2-1),$$

and (5.1) implies that $\varepsilon(x) \in \mathbb{N}\,p_1 + \mathbb{N}\,p_2 \subseteq \sum_{i=1}^{r} \mathbb{N}\,p_i$. $\square$

**Remark 5.3.** Let $P_1, \cdots, P_r$ be the subgroups of $G$ with respectively orders $p_1, \cdots, p_r$. In view of Theorem 5.2, one may wonder if for any $x \in \mathbb{N}G \cap \ker(\varphi)$, there exist $z_i \in \mathbb{Z}G$ with $\varepsilon(z_i) \geq 0$ such that

$$x = z_1\,\sigma(P_1) + \cdots + z_r\,\sigma(P_r). \tag{5.4}$$

This would, of course, imply Theorem 5.2 directly by taking augmentation. Unfortunately, such a representation is not possible in general, if $r \geq 3$. To construct a counterexample, let $|G| = 30$ with $p_1 = 2$, $p_2 = 3$, $p_3 = 5$, and $P_1 = \langle t \rangle$, $P_2 = \langle h \rangle$, $P_3 = \langle g \rangle$. Then, as we saw in (2.5), the element

$$x = t(h+h^2) + g + g^2 + g^3 + g^4 \tag{5.5}$$

lies in $\mathbb{N}G \cap \ker(\varphi)$. Suppose $z_1, z_2, z_3$ exist as in (5.4), with all $\varepsilon(z_i) \geq 0$. Then, from

$$6 = \varepsilon(x) = 2\,\varepsilon(z_1) + 3\,\varepsilon(z_2) + 5\,\varepsilon(z_3),$$

we see that $z_2 = z_3 = 0$ or $z_1 = z_3 = 0$. In the former case, $x = z_1 \cdot (1+t)$; but then $xt = z_1 \cdot (1+t)t = z_1 \cdot (1+t) = x$, which is impossible. Similarly, we see that the case $x = z_2 \cdot (1+h+h^2)$ is impossible as well.

With the result (5.2), it is easy to compute any weight set $W(m)$. We mention explicitly only the special case of *even* integers $m$, which follows directly from (5.2).

**Corollary 5.6.** *Let $m$ be an even integer. Then $W(m)$ is $2\,\mathbb{N}$ when $m$ is a 2-power, and is $\{0, 2, 4, 6, \cdots, p-1, p, p+1, \cdots\}$ when $m$ is not a 2-power and $p$ is the smallest odd prime dividing $m$.*

## 6. Asymmetric Minimal Elements of Smallest Support

We return now to the Lower Bound Theorem 4.8 to give more precise information on the asymmetric minimal elements of the smallest support (resp. weight) in $\mathbb{N}G \cap \ker(\varphi)$. *The notations used in §4 will therefore remain in force throughout this section.*

First we note that the lower bound $p_1(p_2-1) + (p_3-p_2)$ in (4.8) can be written in the more symmetrical form $(p_1-1)(p_2-1) + (p_3-1)$. In the situation of (4.8), this lower bound is the best possible. Indeed, we have seen earlier in (2.5) that there is an asymmetric minimal vanishing sum of distinct $m$th roots of unity, of weight exactly $(p_1-1)(p_2-1)+(p_3-1)$. Transcribing (2.6) in group ring notations, the corresponding asymmetric minimal element in $\mathbb{N}G \cap \ker(\varphi)$ is

$$x(G) := \sigma(P_1^*)\sigma(P_2^*) + \sigma(P_3^*), \tag{6.1}$$

where $P_i$ denotes the subgroup of order $p_i$ in $G$, and $P_i^* := P_i \setminus \{1\}$. Since $x(G)$ realizes the lower bound in (4.8), one wonders naturally about its *uniqueness* (up to similarity). In the following, we shall prove this uniqueness property of $x(G)$. To this end, we must go back to the work in §4, and find out exactly when some of the inequalities there can hold as equalities. We begin with (4.1).

**Proposition 6.2.** *Keep the notations in* (4.1), *and let* $x,\, y \in \mathbb{N}G$ *be such that* $\varphi(x) = \varphi(y)$, $\varepsilon_0(x) \le p_1 - 1$, *and* $y \not\geq x$. *Then:* $\varepsilon_0(y) \ge (p_1 - \varepsilon_0(x))(p_2-1)$, *with equality iff, after a rotation,* $x,\, y$ *are as in* (4.2), *i.e.* $x = c\cdot\sigma(X)$ *and* $y = c\cdot\sigma(X')\,\sigma(P_2^*)$, *where* $c$ *is a positive integer, and* $X,\, X'$ *are two (nonempty) sets forming a partition of* $P_1$.

**Proof.** Applying (4.1), we get the inequality $\varepsilon_0(y) \ge (p_1-\varepsilon_0(x))(p_2-1)$. If $x$ and $y$ have the form described above, then equality holds as we have checked in (4.2). Conversely, assume that $\varepsilon_0(y) = (p_1-\varepsilon_0(x))(p_2-1)$. To pin down the structure of $x$ and $y$, we retrace the steps taken in the proof of (4.1). In particular, we shall use all the notations introduced in that proof, inducting on $r \ge 2$.

*Case 1.* $\varepsilon_0(x_k)+\varepsilon_0(y_j) \ge p_1$. Since $\varepsilon_0(y)$ is exactly $(p_1-\varepsilon_0(x))(p_2-1)$, the various inequalities used in the earlier analysis of this case must all be equalities. This yields very specific information about $x$ and $y$. To begin with, since $\varepsilon_0(y_k)$ was discarded in the earlier estimate of $\varepsilon_0(y)$, we must have $\varepsilon_0(y_k) = 0$, that is, $y_k = 0$. After a rotation (by a power of $g$), we may assume that $k = 0$. For all the other inequalities to be equalities, we must have $p = p_2$, $t = 1$, and (now that $k = 0$) also

$$\varepsilon_0(y_1) = \cdots = \varepsilon_0(y_{p-1}) = p_1 - \varepsilon_0(x_0). \tag{6.3}$$

On the other hand, (4.3) now amounts to $-\varphi(x_0) = \varphi(y_1) = \cdots = \varphi(y_{p-1})$. Therefore, by (3.3)(1), $y_i + x_0 = c_i\cdot\sigma(P_1)$ for each $i \ge 1$, where $c_i \in \mathbb{N}\setminus\{0\}$. Let $X \subseteq P_1$ be the support of $x_0$, and $X' = P_1 \setminus X$. In view of (6.3) and $y_i + x_0 = c_i\cdot\sigma(P_1)$, we see easily that $y_i = c_i\cdot\sigma(X')$ and $x_0 = c_i\cdot\sigma(X)$ for all $i \ge 1$; in particular, all $c_i$'s are equal, say to $c$. We have therefore $x = x_0 = c\cdot\sigma(X)$, and

$$y = y_1 g + \cdots + y_{p-1}g^{p-1} = c\cdot\sigma(X')\,(g+\cdots+g^{p-1}) = c\cdot\sigma(X')\,\sigma(P_2^*).$$

*Case 2.* $\varepsilon_0(x_k)+\varepsilon_0(y_j) \le p_1 - 1$. If $r = 2$, the earlier argument gives $y \ge x$, which is not the case. Hence, we must have $r \ge 3$, and, following through the earlier proof, there exists $\ell$ such that (4.6) holds. Since $\varepsilon_0(y) = (p_1-\varepsilon_0(x))(p_2-1)$, the work on inequalities following (4.6) shows that we must have $\varepsilon_0(x_\ell) = \varepsilon_0(x)$, $\varepsilon_0(y_i) = 0$ for all $i \ne \ell$. This implies that $x = x_\ell g^\ell$ and $y = y_\ell g^\ell$, so after a rotation we may assume that $x,\, y \in \mathbb{N}H$, and we are done by invoking the inductive hypothesis. □

**Corollary 6.4.** *In the notations of* (4.1), *let* $y \in \mathbb{N}G$ *be such that* $\varphi(y) = c_1$ *(a positive integer),* $\varepsilon_0(y) = (p_1 - 1)(p_2 - 1)$, *and* $y \not\geq 1$. *Then* $y = c_1 \cdot \sigma(P_1^*)\,\sigma(P_2^*)$.

**Proof.** For $x := c_1 > 0$, we have $\varphi(x) = \varphi(y)$, $\varepsilon_0(x) = 1 \leq p_1 - 1$, and $y \not\geq x$. By (6.2) applied to $x$, $y$, there exist $u \in G$ and an integer $c > 0$ such that $u \cdot x = c \cdot \sigma(X)$ and $u \cdot y = c \cdot \sigma(X')\,\sigma(P_2^*)$, where $X$, $X'$ are as in (6.2). Since $x = c_1 \in \mathbb{N}$, we must have $c = c_1$ and $X = \{u\}$, so now

$$y = cu^{-1} \cdot \sigma(X')\,\sigma(P_2^*) = c_1 \cdot \sigma(P_1^*)\,\sigma(P_2^*). \quad \square$$

We come now to the main result of this section, which ascertains the uniqueness of asymmetric minimal elements of the smallest support (resp. weight), for a given cyclic group $G$.

**Uniqueness Theorem 6.5.** *Let* $r \geq 3$ *in the notation of* (4.8), *and let* $x$ *be any asymmetric minimal element in* $\mathbb{N}G \cap \ker(\varphi)$. *If either* $\varepsilon(x)$ *or* $\varepsilon_0(x)$ *is equal to* $(p_1 - 1)(p_2 - 1) + (p_3 - 1)$, *then* $x$ *is similar to the element* $x(G)$ *defined in* (6.1).

**Proof.** Since $\varepsilon(x) \geq \varepsilon_0(x) \geq (p_1 - 1)(p_2 - 1) + (p_3 - 1)$ by (4.8), it is sufficient to treat the case $\varepsilon_0(x) = (p_1 - 1)(p_2 - 1) + (p_3 - 1)$. We refer therefore to the proof of (4.8), and retrace the case distinctions there in order to determine the structure of $x$. Since we have at least one strict inequality in *Case 1*, this case cannot occur. In *Case 2*, we are reduced from $\mathbb{N}G \cap \ker(\varphi)$ to $\mathbb{N}H \cap \ker(\varphi)$, so we are done by induction. Thus, it remains only to treat *Case 3*, in which $1 \leq \varepsilon_0(x_i) \leq p_1 - 1$. Here, *Subcase 1* cannot occur since $x$ is asymmetric. Therefore, there must exist an index $j$ as in *Subcase 2*. Looking over the inequality work in that subcase, we see that, for $\varepsilon_0(x) = (p_1 - 1)(p_2 - 1) + (p_3 - 1)$ to hold, we must have $r = 3$, $\varepsilon_0(x_i) = 1$,

$$\varepsilon_0(x_j) \geq (p_1 - \varepsilon_0(x_i))(p_2 - 1) = (p_1 - 1)(p_2 - 1),$$

and also $\varepsilon_0(x_k) = \varepsilon_0(x_i) = 1$ for all $k \neq j$. After a rotation (by a power of $g$), we may assume that $j = 0$. Thus, we have now $x_k = c_k h_k$ for all $k \geq 1$, where $h_k \in H$, and $c_k \in \mathbb{N} \setminus \{0\}$. Since $\varphi(x_0) = \varphi(x_1) = \cdots = \varphi(x_{p-1})$ and $\varphi$ is injective on $H$, we see easily that $c_1 = \cdots = c_{p-1}$, and $h_1 = \cdots = h_{p-1}$. After another rotation (by $h_1^{-1}$), we may therefore assume that $h_1 = \cdots = h_{p-1} = 1$, so $x$ has now the form $x_0 + c_1(g + g^2 + \cdots + g^{p-1})$, with $\varphi(x_0) = c_1$. Clearly, $x_0$ cannot have the identity element $1$ in its support, for otherwise $x \geq 1 + g + \cdots + g^{p-1}$. Also, since $r = 3$,

$$\varepsilon_0(x_0) = \varepsilon_0(x) - (p - 1) = (p_1 - 1)(p_2 - 1) + (p_3 - 1) - (p - 1) = (p_1 - 1)(p_2 - 1).$$

Therefore, (6.4) shows that $x_0 = c_1\,\sigma(P_1^*)\,\sigma(P_2^*)$, and we have

$$x = c_1\,\sigma(P_1^*)\,\sigma(P_2^*) + c_1(g + \cdots + g^{p-1}) = c_1 \cdot x(G).$$

By the minimality of $x$, $c_1$ must be $1$, so $x = x(G)$, as desired. $\square$

The techniques used above for proving the uniqueness of the asymmetric minimal elements of the smallest weight can also be used to analyze asymmetric minimal elements of slightly higher weight. For instance, for weight one higher than the canonical lower bound $(p_1-1)(p_2-1)+(p_3-1)$, one can prove:

**Proposition 6.6.** *Let* $x \in \mathbb{N}G \cap \ker(\varphi)$ *be an asymmetric minimal element of weight* $(p_1-1)(p_2-1)+p_3$. *Then we must have* $p_1 = 2$, $p_2 = 3$, *and* $x$ *is similar to* $t(h+h^2)(1+d)+d^2+d^3+\cdots+d^{p_3-1}$, *where* $t$, $h$, $d$ are suitable generators of the cyclic groups $P_1$, $P_2$ *and* $P_3$. *In particular, such* $x$ *cannot exist if* $|G|$ *is not divisible by* 6.

We shall not go into the details of the proof of this Proposition here. Instead, we offer an explicit illustrative example below.

**Example 6.7.** Let $G$ be a cyclic group of order $m = 30$ as in (5.3), and use the notations there. Then $G = \langle z \rangle$ where $z = thg$, and we have $t = z^{15}$, $h = z^{10}$, $g = z^6$. The map $\varphi : \mathbb{Z}G \to \mathbb{Z}[\zeta_{30}]$ is defined by $\varphi(z) = -\alpha$, where $\alpha := \zeta_{15}$. According to (6.5), the asymmetric minimal element in $\mathbb{N}G \cap \ker(\varphi)$ of the smallest weight is (up to a rotation):

$$x(G) = t(h+h^2)+g+g^2+g^3+g^4 = z^5+z^6+z^{12}+z^{18}+z^{24}+z^{25},$$

with weight $(p_1-1)(p_2-1)+p_3-1 = 6$. Now form the element $x$ of weight 7 in (6.6), using $d = g^2$ as generator for $P_3$ (*cf.* [Ma: p.114]):

$$\begin{aligned} x &= t(h+h^2)(1+d)+d^2+d^3+d^4 \\ &= t(h+h^2)(1+g^2)+g^4+g+g^3 \\ &= z^{25}+z^5+z^7+z^{17}+z^{24}+z^6+z^{18}. \end{aligned}$$

This $x$ corresponds to the vanishing sum

$$\begin{aligned} 0 = \varphi(x) &= -\alpha^{25}-\alpha^5-\alpha^7-\alpha^{17}+\alpha^{24}+\alpha^6+\alpha^{18} \\ &= -\alpha^{10}-\alpha^5-\alpha^7-\alpha^2+\alpha^9+\alpha^6+\alpha^3 \\ &= -\alpha^2(\alpha^8+\alpha^3+\alpha^5+1-\alpha^7-\alpha^4-\alpha). \end{aligned}$$

Since $\alpha$ has degree 8 over $\mathbb{Q}$, the vanishing sum of 30 th roots of unity in parentheses above is clearly minimal. (Incidentally, this shows that $\Phi_{15}(X) = X^8-X^7+X^5-X^4+X^3-X+1$.) Therefore, $x$ is indeed an asymmetric minimal element in $\mathbb{N}G \cap \ker(\varphi)$, arising essentially from the cyclotomic relation satisfied by $\alpha = \zeta_{15}$ over the rationals.

## 7. An Application to Representation Theory

To close this paper, we note the following interesting application of the Main Theorem to the theory of characters of finite groups.

**Theorem 7.1.** *Let $\chi$ be the character of a representation of a finite group $G$ over a field $F$ of characteristic $0$. Let $g \in G$ be an element of order $m = p_1^{a_1} \cdots p_r^{a_r}$ (where $p_1 < p_2 < \cdots$) such that $\chi(g) \in \mathbb{Z}$, and let $t := \chi(1) + |\chi(g)|$. If $\chi(g) \le 0$, then $t \in \sum \mathbb{N}\, p_i$. If $\chi(g) > 0$ and $t$ is odd, then $t \ge \ell$ where $\ell$ ($= p_1$ or $p_2$) is the smallest odd prime dividing $m$.*

**Proof.** Let $D : G \longrightarrow \mathrm{GL}_n(F)$ be the representation in question. Let $A = D(g)$, and let $\alpha_1, \cdots, \alpha_n$ be the eigenvalues of $A$. Then $\alpha_i^m = 1$ for each $i$, and $\chi(g) = \alpha_1 + \cdots + \alpha_n$. Now suppose $s := \chi(g) \in \mathbb{Z}$.

*Case 1.* $s \le 0$. In this case, $\alpha_1 + \cdots + \alpha_n + (-s) \cdot 1 = 0$ is a vanishing sum of $m$th roots of unity, of weight $n - s = \chi(1) + |\chi(g)| = t$. By the Main Theorem (5.2), $t \in \sum \mathbb{N}\, p_i$.

*Case 2.* $s > 0$. In this case, $\alpha_1 + \cdots + \alpha_n + s \cdot (-1) = 0$ is a vanishing sum of $2m$th roots of unity, of weight $n + s = t$. Again by the Main Theorem (5.2), $t \in 2\,\mathbb{Z} + \sum \mathbb{N}\, p_i$. If $t$ is odd, then $m$ must have an odd prime divisor, and if $\ell$ is the smallest odd prime divisor of $m$, then $t \in 2\,\mathbb{Z} + \sum \mathbb{N}\, p_i$ implies that $t \ge \ell$. □

**Examples 7.2.** (1) Let $G = S_8$ and let $\chi$ be the unique irreducible character of degree $7$ on $G$. For $g = (23)(45678) \in G$, we have $\chi(g) = 0$. Here $t = \chi(1) + |\chi(g)| = 7$, which is indeed an $\mathbb{N}$-linear combination of $2$ and $5$ (the prime divisors of the order of $g$). Similarly, if $g' = (123)(45678)$, then $\chi(g') = -1$, and $t' = \chi(1) + |\chi(g')| = 8$, which is an $\mathbb{N}$-linear combination of $3$ and $5$ (the prime divisors of the order of $g'$).

(2) Let $G = \mathrm{SL}(2, 7)$ (a group of order 336), and let $\chi$ be one of the two irreducible characters of degree 6 on $G$. It is known that $\chi(g) = 1$ for some element $g \in G$ of order 14 (see [JL: p.406]). Here $t = \chi(1) + |\chi(g)| = 7$ is odd, and is equal to the smallest odd prime divisor of the order of $m$.

There is apparently no analogue of Theorem 7.1 in charactersitic $p$, even for $p'$-elements $g \in G$ with $\chi(g) = 0$. For instance, if $G = \langle g \rangle$ is a cyclic group of order $4$ and $D$ is the 3-dimensional representation $G \longrightarrow \mathrm{GL}_3(\mathbb{F}_5)$ given by $D(g) = \mathrm{diag}(3, 1, 1)$, then $\chi(g) = 3 + 1 + 1 = 0 \in \mathbb{F}_5$, but the dimension of $D$ is not even an $\mathbb{N}$-linear combination of $2$ (prime divisor of the order of $g$) and $5$ (the characteristic of the ground field).

T. Y. Lam, Mathematics Department, University of California, Berkeley, CA, 94720
*E-mail address*: lam@msri.org

K. H. Leung, National University of Singapore, Singapore 0511